\documentclass[a4paper,11pt]{amsart}
\usepackage[latin1]{inputenc}
\usepackage{xcolor}

\usepackage{multicol}
\usepackage{amsfonts, amssymb, amsmath, amsthm, amscd,epsfig,mathrsfs}
\usepackage{latexsym}
\usepackage{graphicx}

\usepackage{cancel}
\usepackage[normalem]{ulem}

\usepackage[all]{xy}

\usepackage{enumerate}

\usepackage{verbatim}

\usepackage{hyperref}

\DeclareMathOperator{\mult}{mult}
\DeclareMathOperator{\ord}{ord}

\DeclareMathOperator{\Sing}{Sing}
\DeclareMathOperator{\Spec}{Spec}

\newcommand{\Diff}{\mathit{Diff}}
\newcommand{\G}{{\mathcal G}}

\renewcommand{\L}{{\mathcal L}}

\newcommand{\p}{{\mathfrak p}}

\newcommand{\LX}{\mathcal{L}(X)}
\newcommand{\LXP}{\mathcal{L}(X')}

\newcommand{\Gn}{\G^{(n)}}

\newtheorem{Thm}{Theorem}[section]      %Numeraci�n dentro de secciones

\theoremstyle{definition}
\newtheorem{Parrafo}[Thm]{\ }
\newtheorem{Def}[Thm]{Definition}  
\newtheorem{Def-Prop}[Thm]{Definition-Proposition}

\theoremstyle{remark}
\newtheorem{Rem}[Thm]{Remark}

\newtheorem{Ex}[Thm]{Example}%[chapter]

\numberwithin{equation}{Thm}

\title{Finite morphisms and Nash multiplicity sequences}
\author{A. Bravo, S. Encinas}
\thanks{The authors were partially supported by PGC2018-095392-B-I00;
The first author was partially supported from the Spanish Ministry of Economy and Competitiveness, through the ``Severo Ochoa'' Programme for Centres of Excellence in R\&D (SEV-2015-0554)}

\keywords{Rees algebras. Resolution of Singularities. Arc Spaces. Finite morphisms.}
\subjclass[2010]{14E15, 14E18}

\AtEndDocument{\bigskip{\footnotesize
		\textsc{Depto. Matem\'aticas,
			Facultad de Ciencias, Universidad Aut\'onoma de Madrid 
			and Instituto de Ciencias Matem\'aticas CSIC-UAM-UC3M-UCM, Canto Blanco 28049 Madrid, Spain} \par 
		\textit{E-mail address}, A. Bravo: \texttt{ana.bravo@uam.es} \par
		\addvspace{\medskipamount}
		\textsc{Depto. Matem\'atica Aplicada,
			and IMUVA, Instituto de Matem\'aticas.
			Universidad de Valladolid} \par
		\textit{E-mail address}, S. Encinas: \texttt{santiago.encinas@uva.es}
}}

\begin{document}

\begin{abstract}
We study finite morphisms of varieties and the link between their top multiplicity loci under certain assumptions.
More precisely, we focus on how to determine that link in terms of the spaces of arcs of the varieties. 
\end{abstract}

\maketitle

%\tableofcontents

\section{Introduction}

%\noindent 
The multiplicity of a variety $X$ at a singular point can be understood as a measure of the singularity: $X$ is regular if and only if the multiplicity at any of its points is 1.   Observe also that    the multiplicity defines an upper-semi continuous function on $X$. As a consequence, if $m$ is the maximum multiplicity of $X$, then the set of points with multiplicity $m$, $F_m(X)$, is closed. 

\

The multiplicity does not increase when blowing up along regular centers contained in $F_m(X)$ (see \cite{Orbanz} or \cite{Dade}). 
Motivated by this fact, we say that a closed regular subscheme $Y\subset X$ is {\em $F_m$-permissible} if $Y \subset F_m(X)$. A blow up at an $F_m$-permissible centers be called a \emph{$F_m$-permissible blow up}. 

\

We will say that a sequence of $F_m$-permissible blow ups 
\begin{equation}
\label{simply_1}
\xymatrix @R=0pt @C=30pt {
	X=X_0 & \ar[l]_-{\pi_1} X_1 & \ar[l]_-{\pi_2} \ldots & \ar[l]_-{\pi_l} X_l,
} 
\end{equation}
such that
\begin{equation}
\label{simply_2}
\max\mult(X_0)
= \dotsb = \max\mult(X_{l-1}) > \max\mult(X_l)
\end{equation}
is a \emph{simplification of the multiplicity, $F_m(X)$}.

\

 When  $X$ is defined over a field of characteristic zero, one can achieve  a resolution of singularities of $X$   by iterating succesive simplifications of the multiplity of $X$    (cf. \cite{V}). This is shown by attaching a suitably defined Rees algebra $\mathcal{G}_X$ to the closed set $F_m(X)$. This Rees algebra helps to describe the closed set  $F_m(X)$ and in addition provides sufficient information to define  invariants that ultimately lead to the construction of    a sequence like (\ref{simply_1}) such that (\ref{simply_2}) holds.   More  details about this fact will be given in the following paragraphs
 and forthcoming sections.

\ 

In this paper we are interested in the study of a class of finite morphisms between varieties and the link between their top multiplicity loci. More precisely, let $k$ be a perfect field    
and let  $\beta: X'\to X$ be a finite (dominant) morphism of $k$-varieties of generic rank $r$.  Suppose that the maximum multiplicity at the points of $X$ is $m$. Then the maximum multiplicity at the points of $X'$ is bounded above by $rm$. When this upper bound is attained we say that $\beta: X'\to X$ is {\em transversal}. When $\beta$ is transversal there is an interesting link between the (closed) set of points of multiplicity $rm$ in $X'$, $F_{rm}(X')$, and the top multiplicity locus of $X$, $F_m(X)$. For instance, it can be proven that $F_{rm}(X')$ is homeomorphic to $\beta(F_{rm}(X'))$, which, in addition, sits inside $F_m(X)$. This and other properties of transversal morphisms have been studied in the context of constructive resolution of singularities in \cite{COA}. 

\

The algebra $\mathcal{G}_X$ from above can always  be defined for varieties over perfect fields \cite{V07}.
When $\beta:X'\to X$ is transversal  it can be shown that there is  an extension of the Rees algebras $\mathcal{G}_X\subset\mathcal{G}_{X'}$ associated to $F_m(X)$ and $F_{rm}(X')$ respectively (cf. \cite{Abad}). 

\

Suppose now that $F_{rm}(X')$ is homeomorphic to $F_m(X)$.
A natural question is to wonder  if there is a    link between the simplifications of the multiplicities for $F_{rm}(X')$ and $F_m(X)$.

\

 In the characteristic zero case, one of the results of \cite{COA} says that if the extension
$\mathcal{G}_X\subset\mathcal{G}_{X'}$ is finite then a simplification of $F_m(X)$ induces naturally a simplification of $F_{rm}(X')$ and vice versa. However in positive characteristic there are examples where the extension 
$\mathcal{G}_X\subset\mathcal{G}_{X'}$ is finite but there is not such a  strong link between $F_m(X)$ and $F_{rm}(X')$   (see \cite[Example 7.5]{COA}).  

\

 Our purpose is to study transversal morphisms using properties of the arc spaces of both $X$ and $X'$, when they are varieties defined over a perfect field.
More precisely, we will be looking at the Nash mutiplicity sequences of the arcs with center in the top multiplicity loci of the varieties. The main result is Theorem \ref{MainTheorem},
  where we give an equivalent condition to the finiteness of the extension $\mathcal{G}_X\subset\mathcal{G}_{X'}$ in terms of the spaces of arcs of $X$ and $X'$.

\

In the following lines we give more details, definitions, the motivation for our problem and its relation with constructive resolution of singularities.

\

\noindent {\bf Constructive resolution of singularities and multiplicity}
\medskip

\noindent After Hironaka's Theorem on resolution of singularities in characteristic zero \cite{Hir}, a series of algorithms of resolution were found (\cite{B-M}, \cite{V1}, and \cite{V2};
see also \cite{Br_E_V}, \cite{E_V97} and \cite{E_Hau}). An algorithmic resolution of singularities consists on describing a procedure to construct, step by step, a sequence of blow ups that leads to the resolution of a given variety $X$, 
\begin{equation}
\label{constructive}
X=X_0\leftarrow X_1 \leftarrow \ldots \leftarrow X_n=T.
\end{equation}

Roughly speaking, to find a sequence like (\ref{constructive}) one uses the so called {\em resolution functions defined on varieties}. These are upper semi-continuous functions, 
$$ \begin{array}{rrcl} 
f_X: & X & \to & (\Lambda, \geq) \\
& \xi & \mapsto & f_X(\xi),
\end{array}$$
that are constant if and only if the variety is regular and whose maximum value, $\max f_X$, achieved in a closed regular subset $\underline{\text{Max}} f_X$, selects the center to blow up. Thus the sequence (\ref{constructive}) is defined so that 
$$\max f_{X_0}> \max f_{X_1} > \ldots > \max f_{X_n},$$
where $\max f_{X_i}$ denotes the maximum value of $f_{X_i}$ for $i=0,1,\ldots, n$.

\

Usually, $f_X$ is defined, at each point, as a sequence of rational numbers, the first set of coordinates being the Hilbert-Samuel function at the point (see \cite{Br_E_V}) or the multiplicity (see \cite{V}). For the purposes of this paper we will be paying attention to the later.  
  Therefore we will be considering  a  resolution function on $X$  as the following: 
\begin{equation}
\label{resol_coord}
f_X(\xi)=(\text{mult}_{X}(\xi),\ldots). 
\end{equation}
And we will be achieving a desingularization of $X$ by concatenating successive simplifications of the multiplicity of $X$. 

\

\noindent {\bf On refinements of the multiplicity}
\medskip

\noindent   Now let us say a word about the other coordinates of $f_X$ in (\ref{resol_coord}). Even though the multiplicity is an upper-semi continuous function on $X$, it usually does not define a resolution function. For instance the closed set $F_m(X)$ may not be even regular.  Therefore, in order to construct a resolution function we need to find refinements of the multiplicity. 
 These are defined by using {\em local presentations of the multiplicity} (this was studied in 	\cite[\S 7.1]{V}).

\

Roughly speaking  by a local presentation of the multiplicity in a neighborhood of a point $\xi\in F_m(X)$ we mean   that locally, in an \'etale neighborhood of a point $\xi$, which we denote again by $X$ for simplicity, one can find an embedding of $X$ in some smooth scheme, $V$, together with a set of weighted equations that (locally) describe $F_m(X)$ (see Example \ref{Ex:HiperRees} for the case of a hypersurface).
The information given by such (finite) set of weighted equations is expressed in terms of a Rees algebra $\G$ defined on $V$  (see section~\ref{Rees_Algebras}), and we refer to the pair $(V,\G)$ as a {\em local presentation of $F_m(X)$}.

\

Local presentations are not unique, i.e., there may be different embeddings and different Rees algebras that provide local presentations of $F_m(X)$. However, it   can be proven that they all lead to the same resolution function \cite[Theorem 26.5]{Br_V2}.

\

In addition, it can be shown that the restriction of $\G$ to $X$, $\G_X$, is unique up to integral closure, (cf. \cite{Abad}). We will say that $\G_X$ is {\em the ${\mathcal O}_X$-Rees algebra attached to $F_m(X)$} in a neighborhood of $\xi$. We refer to section \ref{Rees_Algebras} for precise definitions and 
statements regarding Rees algebras and local presentations.

\

 When the characteristic is zero, the pair $(V, \G)$ provides all the information needed to construct a simplification of $F_m(X)$ locally in a neighborhood of $\xi\in F_m(X)$; in other words, the remaining coordinates of $f_X$ at the points in $F_m(X)$ are determined by $(V, \G)$ (see (\ref{resol_coord})). For instance, if $X$ is a $d$-dimensional variety, then 
\begin{equation}
\label{resol_coord_1}
f_X(\xi)=(\text{mult}_{X}(\xi), \ord^{(d)}_{X}(\xi), \ldots), 
\end{equation}
where $\ord^{(d)}_{X}(\xi)$ is a rational number that we refer to as {\em Hironaka's order function in dimension $d$}. This number is obtained by performig some sort of elimination of variables on $(V, \G)$ (see \cite{Abad-PhD}, \cite{Br_V}, \cite{V07}), and it can be seen as a refinement of the multiplicity that leads to the construction of a resolution function.

\

\noindent{\bf Nash multiplicity sequences and constructive resolution}
\medskip

\noindent The rational number $\ord_X^{(d)}(\xi)$ from above can be defined whenever $k$ is a perfect field (in positive characteristic it does not provide enough information to define a resolution function).   See \cite{Abhy1}, \cite{Abhy2}, \cite{Benito_V}, \cite{Cos_Pilt1,Cos_Pilt2}, \cite{Cut}, \cite{K_M1}, \cite{Lipman3} for results on resolution in positive characteristic. 

\

 In \cite{Br_E_P-E} and \cite{BEP2}  with B. Pascual-Escudero, we showed that $\ord_X^{(d)}(\xi)$, can be read from the set of arcs with center $\xi$, $\L(X, {\xi})$. 
 To this end we worked with the so called {\em Nash multiplicity sequences of arcs}. These were introduced by M. Lejeune-Jalabert in \cite{L-J} for the case of a germ of a point of a hypersurface, and generalized afterwards by H. Hickel in \cite{Hickel93} and  \cite{Hickel05}.

\

Given a point $\xi\in \Sing(X)$ of multiplicity $m>1$, and an arc $\varphi \in \L(X, \xi)$ the {\em sequence of Nash multiplicities of $\varphi$} is a non-incressing sequence of integers, 
\begin{equation}
\label{primera_Nash}
m=m_0\geq m_1\geq m_2 \geq \ldots
\end{equation}
where $m_0=m$ is the multiplicity at the point $\xi$, and the rest of the numbers in the sequence can be interpreted as a {\em refinement of the ordinary multiplicity at $\xi$ along the arc $\varphi$}.
 See section~\ref{Jet_Arcs}, specially diagram (\ref{intro:diag:Nms})   for details on the definition of this sequence. 

\

Suppose that $\varphi$ is a $K$-arc, with $K\supset k$, which defines a morphism 
$\varphi:{\mathcal O}_{X,\xi}\to K[[t]]$. 
When the generic point of $\varphi$ is not contained in the stratum of multiplicity $m$ of $X$, then there is some subindex $l\ge 1$ in sequence (\ref{primera_Nash}) for which $m_l<m_0$. We are interested in the first subindex for which the inequality holds and call it {\em the persistance of the arc $\varphi$}, $\rho_{X,\varphi}$. To eliminate the impact of the order of the arc at the point, we normalize the persistance setting 
\begin{equation}
\label{rho_norm_1}
\overline{\rho}_{X,\varphi}=\frac{\rho_{X,\varphi}}{\nu_t(\varphi)},
\end{equation}
where $\nu_t(\varphi)$ denotes the order of the image by $\varphi$, of the defining ideal of $\xi$ at the regular local ring $K[[t]]$. We work simultaneously with another invariant which is refinement of the persistance: the {\em ${\mathbb Q}$-persistance}, which we denote by $r_{X,\varphi}$, and its normalized version $\overline{r}_{X,\varphi}$. In fact, the two invariants are related since for a given arc $\varphi$ it can be shown that 
\begin{equation} 
\label{rho_r} 
\rho_{X,\varphi}= \lfloor r_{X,\varphi}\rfloor \ \text{ and } \ r_{X,\varphi} = \lim _{n\rightarrow \infty }\frac{\rho _{X,\varphi _{n}}}{n} \in {\mathbb Q}_{\geq 1},
\end{equation}
where for each $n\geq 1$, $\varphi _n=\varphi \circ i_n$ and $i_n^*:K[[t]]\longrightarrow K[[t]]$
  is the $K$-morphism mapping $t$ to $t^n$. 
\medskip

Using these definitions, in joint works with B. Pascual-Escudero we showed that Hironaka's order function in dimension $d$ (\ref{resol_coord_1}) can be read from the Nash multiplicity sequences of the arcs in $X$: 

\medskip

\begin{Thm}\label{principal} \cite[Theorem 3.6]{Br_E_P-E}, \cite[Theorem 6.1]{BEP2} 
	Let $X$ be a $d$-dimensional algebraic variety defined over a perfect field $k$, and let $\xi \in F_m(X)$. Then 
	\begin{equation}\label{desigualdad_i}
	\ord_X^{(d)}(\xi)\leq \inf_{\varphi\in {\L(X, {\xi})}}\{\overline{r}_{X,\varphi}\}=\inf_{\varphi\in \L(X,{\xi})}\left\{\frac{1}{\nu_t(\varphi)}\lim_{n\to \infty}\frac{\rho_{X,\varphi_{n}}}{n}\right\}.
	\end{equation}
	Moreover, the infimum is a minimum, i.e., there is some arc $\eta\in \L(X,{\xi})$ such that: 
	\begin{equation}
	\label{igualdad_divisorial_i}
	\ord_X^{(d)}(\xi)=\overline{r}_{X,\eta}=\frac{1}{\nu_t(\eta)}\lim_{n\to \infty}\frac{\rho_{X,\eta_{n}}}{n}.
	\end{equation}
\end{Thm}

Note that for the definition of $\ord^{(d)}_X(\xi)$ (\ref{Th:PresFinita}), it is necessary to find a suitable \'etale neighborhood of $\xi$, a local embedding in a smooth scheme,
and the construction of a convenient Rees algebra. 
A consequence of Theorem~\ref{principal} is that $\ord^{(d)}_X(\xi)$ can be defined without using \'etale topology   and only studying properties of its space of arcs.  
Moreover the arc $\eta$ realizing the minimum in (\ref{igualdad_divisorial_i}) can be choosen, and constructed explicitly,  being   fat and divisorial \cite[Theorem~6.3]{BEP3}.   In other words, the refinement of the multiplicity for the resolution function in (\ref{resol_coord}) can be obtained by studying sequences of Nash multipicities sequences in ${\mathcal L}(X)$. 

\

\noindent {\bf Simplifications of the multiplicity and finite morphisms between singular varieties} 
\medskip

\noindent In \cite{COA} the following problem is studied. Let $\beta : X' \to X$ be a finite and dominant morphism of singular algebraic varieties over a perfect field $k$, and suppose, as before, that the maximum multiplicity of $X$ is $m$. Assume that $r = [K(X'):K(X)]$ is the generic rank of $\beta$. Then by Zariski's multiplicity formula for finite projections it follows that the maximum multiplicity on $X'$ is bounded above by $rm$.  When $F_{rm}(X')\neq \emptyset$, we say that the morphism $\beta : X' \to X$ is \emph{transversal}. In such case  $F_{rm}(X')$ is mapped homeomorphically to a closed subset in   $F_m(X)$.  

\

Transversality is preserved by $F_{rm}$-permissible blow ups. In other words, any $F_{rm}$-permissible blow up on $X'$, $X'\leftarrow X'_1$ induces an $F_{m}$-permissible blow up on $X$, $X\leftarrow X_1$, and there exists a finite dominant morphism $\beta_1:X'_1\to X_1$ making the diagram
\begin{equation*}
\xymatrix@R=15pt{
			X' \ar[d]^\beta &
			X'_1 \ar[l] \ar[d]^{\beta_1} \\
			X &
			X_1, \ar[l]
} 
\end{equation*}
commutative; if in addition $F_{rm}(X'_1)\neq\emptyset$ then the morphism $\beta_1$ is transversal (see \cite[Theorem 4.4]{COA}). 
We refer to section \ref{transversalidad} for more details and precise statements.

\

It is natural to study conditions under which, given a transversal morphism $\beta : X' \to X$, the set $F_{rm}(X')$ is mapped surjectively onto $F_m(X)$, in such a way that $F_{rm}(X')$ and $F_m(X)$ are homeomorphic and, in addition, the condition is preserved by sequences of $F_m$-permissible blow ups. In \cite{COA} these morphisms are called \emph{strongly transversal}  (see Definition~\ref{def:strongly-transversal}). If $\beta : X' \to X$ is strongly transversal then a simplification of the multiplicity of $X'$ induces a simplification of the multiplicity of $X$ and vice versa.

\

\noindent {\bf Strong transversality, Rees algebras and Nash multiplicity sequences} 
\medskip

\noindent As indicated above, there are Rees algebras attached to the maximum multiplicity loci of $X$ and $X'$, $\G_X$ and $\G_{X'}$. When the characteristic is zero, strong transversality can be characterized in terms of the algebras $\G_X$ and $\G_{X'}$. In fact, the following theorem holds over perfect fields: 

\begin{Thm}
\cite[Theorem 7.2]{COA}
Let $\beta : X' \to X$ be a transversal morphism of generic rank $r$ between two singular algebraic varieties defined over a perfect field $k$. Then:
\begin{enumerate}
	\item If $\beta : X' \to X$ is strongly transversal then the inclusion $\G_X\subset \G_{X'}$ is finite; 
	
	\item If $k$ is a field of characteristic zero, then the converse holds. Namely, if the inclusion $\G_X \subset \G_{X'}$ is finite, then $\beta : X' \to X$ is strongly transversal.
\end{enumerate}
\end{Thm}

\begin{Rem}\label{contenidos_algebras}
Recall that the Rees algebras $\G_X$ and $\G_{X'}$ are only defined locally in \'etale topology.
However, as we will see in section \ref{proof}, given a point $\xi \in F_m(X)$, one can find an \'etale neighborhood of $X$ at $\xi$, $\widetilde{X} \to X$, where the intrinsic algebra $\G_{\widetilde{X}}$ associated to $\widetilde{X}$, as well as the intrinsic algebra $\G_{\widetilde{X}'}$ associated to $\widetilde{X}' = X' \times_X \widetilde{X}$ are defined.

It is in this setting that there is an inclusion $\G_{\widetilde{X}'} \subset \G_{\widetilde{X}}$, and in which we compare these algebras. See also \cite[Remark 7.3]{COA}.

Also, transversality is preserved by \'etale  base change (\ref{etale_transversal}).
\end{Rem}

Thus, when the characteristic is positive, strong transversality implies that $\G_{X'}$ is integral over $\G_X$ but the converse may fail (see \cite[Example 7.5]{COA} for a counterexample in the latter case). It is natural to ask what piece of information is encoded if the containment $\G_X \subset \G_{X'}$ is finite.

\

The main theorem in this paper says 
that this condition can be expressed in terms of the Nash multiplicity sequences: 

\begin{Thm} \label{MainTheorem} 
Let $\beta : X' \to X$ be a transversal morphism of generic rank $r$ between two singular algebraic varieties defined over a perfect field $k$, and let $\beta_{\infty}: \LXP \to \LX$ be the induced morphism. 

 Let $m$ be the maximum multiplicity of $X$ and assume that $F_{rm}(X')$ is homeomorphic to $F_m(X)$.
Then,  
the inclusion $\G_X\subset \G_{X'}$ is finite if and only if for each arc $\varphi'\in {\mathcal L}(X')$ with center in $F_{rm}(X')$, we have the following equality of persistances:
$$\rho_{X',\varphi'}=\rho_{X,\beta_{\infty}(\varphi')}.$$
\end{Thm}

 In particular, when the characteristic of $k$ is zero, a simplification of the multiplicity of $X'$ induces a simplification of the multiplicity of $X$ and vice versa if and only if for each arc $\varphi'\in {\mathcal L}(X')$, $\rho_{X',\varphi'}=\rho_{X,\beta_{\infty}(\varphi')}$.  On the other hand, when the characteristic is positive, the result  says that the finiteness of the extension  $\G_X \subset \G_{X'}$  indicates a strong link between the  Nash multiplicity sequences of arcs with center  at the top multiplicity loci of both $X$ and $X'$. 

\ 

\medskip

\noindent{\bf On the organization of the paper} 
\medskip

The paper is organized as follows.
In Section~\ref{transversalidad} we review  the meanings of transversality and strong transversality for finite morphisms.
In section~\ref{Rees_Algebras} we recall  the notion of  Rees algebra     and its  use in constructive resolution. In this setting, Rees algebras   are graded algebras which locally,  in \'etale topology, describe the multiplicity maximum locus of a variety (in some strong sense that will be made precise in the section). 
Section \ref{Jet_Arcs}  is devoted to recalling the definition of the Nash multiplicity sequence of an arc, and the concept of persistance associated to an arc in the variety.
Finally   Theorem~\ref{MainTheorem} is proven in  section~\ref{proof}.
 
\medskip

{\em Acknowledgements.} We profited from conversations with C. Abad, A. Benito, B. Pascual-Escudero and O. E. Villamayor.

\section{Transversality and strong transversality} \label{transversalidad} 

As indicated in the introduction, we are interested in studying certain finite morphisms between singular variteties. We will start by recalling Zariski's multiplicity formula for finite projections.
Let $(R,\mathfrak{m})$ be a local Noetherian ring and let $\mathfrak{a}\subset R$ be an $\mathfrak{m}$-primary ideal. We denote by $e_R(\mathfrak{a})$ the multiplicity of $R$ with respect to the ideal $\mathfrak{a}$.
The multiplicity of a Noetherian scheme $X$ at a point $\xi\in X$ is defined as that of the local ring ${\mathcal O}_{X,\xi}$ at its maximal ideal. Zariski's multiplicity formula is stated in the following Theorem:

\begin{Thm}\cite[VIII, Theorem 24, Corollary 1]{Z-SII}\label{MultForm} Let 
	$(A,{\mathfrak m})$ be a local domain and let $C$ be a finite extension of $A$. Let $K$ denote the quotient field of $A$, and let $L = K \otimes_A C$. Let ${\mathfrak n}_1, \ldots, {\mathfrak n}_r$ denote the maximal ideals of the semi-local ring $C$, and assume that $\dim C_{{\mathfrak n}_i}=\dim C$ for $ i=1, \dots ,r$. Then
	\[e_A({\mathfrak m})[L:K] = \sum_{1\leq i \leq r} e_{C_{{\mathfrak n}_i}}({\mathfrak m} C_{{\mathfrak n}_i}) [k_i:k],\]
	where $k_i$ is the residue field of $ C_{{\mathfrak n}_i}$, $k$ is the residue field of $(A,{\mathfrak m})$, and $[L:K]= \dim_K L$.
	
\end{Thm}

Let $X$ be an irreducible algebraic variety over a perfect field $k$, and let $X'$ be an equidimensional algebraic variety over $k$.
Denote by $K$ the field of rational functions of $X$ and let $L$ be the total ring of fractions of $X'$.
If $\beta : X' \to X$ is a finite and dominant $k$-morphism, then by Zariski's formula (\ref{MultForm})
\begin{equation}\label{eq222}
\max\mult(X') \leq [L:K] \cdot \max\mult(X).
\end{equation}

\begin{Def} \label{def:transversal}
\cite[Definition 2.5]{COA}
With the previous notation,
we will say that $\beta : X' \to X$ is \emph{transversal} if:
\begin{equation}\label{eq2222} 
	\max\mult(X')= [L:K] \cdot \max\mult(X).
\end{equation} 
\end{Def}

\begin{Rem} \label{condicion_estrella}
Assume that we are in the affine case, $X=\Spec(B)$ and $X'=\Spec(B')$, where $B$ and $B'$ are $k$-algebras and
the finite morphism $\beta:X'\to X$ is given by a finite extension $B\to B'$. Let $P\in\Spec(B')$ be a point and set $\mathfrak{p}=P\cap B\in\Spec(B)$.
Then the equality:
\begin{equation} \label{IgualMultB}
e_{B'_P}(PB'_P)= e_{B_{\p}}(\p B_{\p})[L:K]
\end{equation}
holds if and only if the following three conditions hold simultaneously:
\begin{enumerate}
	\item[(i)] $P$ is the only prime in $B'$ dominating ${\p}$ (i.e., $B'_P=B'\otimes_BB_{\p}$);	
	\item[(ii)] $B_{\p}/{\p}B_{\p}=B'_P/PB'_P$;
	\item[(iii)] $e_{B'_P}(\p B'_P)=e_{B'_P}(PB'_P)$.
\end{enumerate}
In particular, condition (\ref{IgualMultB}) necessarily holds for all primes $P\subset B'$ with multiplicity $rm$, where $m$ is the maximum multiplicity in $\Spec(B)$, and $r=[L:K]$.

Now, suppose that $B$ and $B'$ are formally equidimensional locally at any prime. Then, condition (iii) is equivalent to saying that $\p B'_P$ is a reduction of $P B'_P$, i.e., that the ideal $P B'_P$ is integral over $\p B'_P$ (cf. \cite{Rees}).
\end{Rem}

\begin{Rem}
	\label{resumen} If $\beta: X'\to X$ is transversal, then it can be shown that: 
	\begin{enumerate} 
		\item $\beta (F_{rm}(X')) \subset F_{m}(X)$;
		\item $F_{rm}(X')$ is homeomorphic to $\beta (F_{rm}(X'))$; 
		\item If $Y \subset F_{rm}(X')$ is an irreducible regular closed subscheme, then $\beta(Y)\subset F_{m}(X)$ is an irreducible regular closed subscheme; 
		\item If $Z \subset F_{m}(X)$ is an irreducible closed regular subscheme, and if $\beta^{-1}(Z)_{\text{red}}\subset F_{rm}(X')$, then $\beta^{-1}(Z)_{\text{red}}$ is regular. 
	\end{enumerate}
See \cite[Proposition 2.7 and Corollary 2.8]{COA}.

\end{Rem}
\begin{Parrafo}\label{local_transformations} {\bf Local transformations.}
We will see that transversality is stable under permissible blow ups and other special morphisms that play an important role in resolution of singularities.

A morphism $X_1 \to X$ is an \emph{$F_m$-local transformation} if it is of one of the following types:
\begin{enumerate}[\textup{(}i\textup{)}]
	\item The blow up of $X$ along a regular center $Y$ contained in $F_m(X)$. This will be called an \emph{$F_m$-permissible blow up}. In this case we will also say that $Y$ is an \textit{$F_m$-permissible center}.
	\item An open restriction, i.e., $X_1$ is an open subscheme of $X$. In order to avoid trivial transformations, we will always require $X_1 \cap F_m(X) \neq \emptyset$.
	\item The multiplication of $X$ by an affine line, $X_1 = X \times \mathbb{A}^1_k$.
\end{enumerate}
Note that, in either case $\max\mult(X) \geq \max\mult(X_1)$. A sequence of transformations, 
\[ \xymatrix {
	X = X_0 &
	X_1 \ar[l]_-{\phi_1} &
	\dots \ar[l]_-{\phi_2} &
	X_N \ar[l]_-{\phi_{N}}
}, \]
is an \textit{$F_m$-local sequence on $X$} if $\phi_i$ is an $F_m$-local transformation of $X_{i-1}$ for $i = 1, \dotsc, N$, and
\begin{equation*}
m = \max\mult(X_0)
= \dotsc = \max\mult(X_{N-1})
\geq \max\mult(X_N) .
\end{equation*}

The morphisms defined in (i)-(iii) above are the starting point for the definition of some of the fundamental invariants in resolution related to the so called {\em Hironaka's trick} \cite[Proposition 7.3]{E_V97}.
\end{Parrafo}

\begin{Thm} \label{trans_permanencia} \cite[Theorem 4.4, Remark 4.5]{COA}.
Let $X$ be an algebraic variety with maximum multiplicity $m$ and let $\beta : X' \to X$ be a transversal morphism of generic rank $r$. Then:
\begin{enumerate}[\textup{(}i\textup{)}]
\item An $F_{rm}$-permissible center on $X'$, $Y' \subset F_{rm}(X')$, induces an $F_m$-permissible center on $X$, $Y = \beta(Y') \subset F_m(X)$,
and a commutative diagram of blow ups of $X$ at $Y$, $X\leftarrow X_1$, and of $X'$ at $Y'$, $X'\leftarrow X_1'$, as follows,
\[ \xymatrix@R=15pt{
			X' \ar[d]^\beta &
			X'_1 \ar[l] \ar[d]^{\beta_1} \\
			X &
			X_1, \ar[l]
} \]
where $\beta_1$ is finite of generic rank $r$. In addition, if $F_{rm}(X'_1) \neq \emptyset$, then $F_m(X_1) \neq \emptyset$, and the morphism $\beta_1$ is transversal.

\item Any $F_{rm}$-local sequence on $X'$, 
$X' \leftarrow X'_1 \leftarrow\dotsb\leftarrow X'_{N-1} \leftarrow X'_N,$
induces $F_{m}$-local sequence on $X$, and a commutative diagram as follows,
\begin{equation*} %\label{eq:Frs-finite-diagram}
\xymatrix@R=15pt{
			X' \ar[d]^{\beta} &
			X'_1 \ar[l] \ar[d]^{\beta_1} &
			\dotsb \ar[l] &
			X'_{N-1} \ar[l] \ar[d]^{\beta_{N-1}} &
			X'_N \ar[l] \ar[d]^{\beta_N} \\
			X &
			X_1 \ar[l] &
			\dotsb \ar[l] &
			X_{N-1} \ar[l] &
			X_N, \ar[l]
}
\end{equation*}
where each $\beta_i$ is finite of generic rank $r$. Moreover, if $F_{rm}(X'_N) \neq \emptyset$, then $F_m(X_N) \neq \emptyset$, and the morphism $\beta_N$ is transversal.
\end{enumerate}
\end{Thm}

\begin{Def} \label{def:strongly-transversal} 
\cite[Definition 4.8]{COA}
We will say that a transversal morphism of generic rank $r$, $\beta : X' \to X$, is \emph{strongly transversal} if $F_{rm}(X')$ is homeomorphic to $F_m(X)$ via $\beta$, and every $F_{rm}$-local sequence over $X'$, $X' \leftarrow X'_1 \leftarrow \dotsb \leftarrow X'_N$,
	induces an $F_m$-local sequence over $X$ and a commutative diagram as follows,
	\begin{equation} \label{diag:def-strong-transversal}
	\xymatrix@R=15pt{
		X' \ar[d]^{\beta} &
		X'_1 \ar[l] \ar[d]^{\beta_1} &
		\dotsb \ar[l] &
		X'_{N-1} \ar[l] \ar[d]^{\beta_{N-1}} &
		X'_N \ar[l] \ar[d]^{\beta_N} \\
		X &
		X_1 \ar[l] &
		\dotsb \ar[l] &
		X_{N-1} \ar[l] &
		X_N, \ar[l]
	} \end{equation}
	where each $\beta_i$ is finite of generic rank $r$ and induces a homeomorphism between $F_{rm}(X'_i)$ and $F_m(X_i)$. In this case we will also say that $F_{rm}(X')$ is \emph{strongly homeomorphic} to $F_m(X)$. Note in particular that this definition yields $F_{rm}(X'_N) = \emptyset$ if and only if $F_m(X_N) = \emptyset$.
\end{Def}

Notice that in Theorem \ref{MainTheorem}  the hypothesis is only that $F_{rm}(X')$ is homeomorphic to $F_m(X)$. This is weaker than saying that the sets $F_{rm}(X')$ and $F_m(X)$ are strongly homeomorphic.  

\

From the point of view of resolution of singularities, $\beta : X' \to X$ is strongly transversal if and only if the processes of simplification of the multiplicity of $X'$ and $X$ are equivalent. Now observe that if $\varphi'\in \L(X')$ has center in $F_{rm}(X')$, then $\varphi:=\beta_{\infty}(\varphi')$ is an arc with center in $F_m(X)$. If we consider the sequence of multiplicities of Nash as a refinement of the usual multiplicity, it is natural to compare the persistances of $\varphi'$ and $\varphi$. 

\begin{Rem}
\label{etale_transversal}
In some of our arguments we will work in \'etale topology, and it is worth noticing that transversality and strong transversality   are  preserved after considering \'etale change of basis.
Suppose we are given a transversal morphism $X'\to X$ and an \'etale morphism $\widetilde{X}\to X$. Then it can be checked that the induced morphism $\widetilde{X}\times_X X'\to \widetilde{X}$ is transversal again.
In the sense that 
equality (\ref{eq2222}) is preserved by base change,
replacing $K$ by the total ring of fractions of $\widetilde{X}$. 
\end{Rem}

\section{Rees algebras}\label{Rees_Algebras}

The stratum defined by the maximum value of the multiplicity function of a variety can be described using equations with  weights (\cite{V}).
The same occurs with the Hilbert-Samuel function (\cite{Hir1}). Along this section we will see that 
Rees algebras are natural objects to work with this setting. 

\begin{Def}
	Let $R$ be a Noetherian ring. A \textit{Rees algebra $\mathcal{G}$ over $R$} is a finitely generated graded $R$-algebra
	$\mathcal{G}=\oplus _{l\in \mathbb{N}}I_{l}W^l\subset R[W]$
	for some ideals $I_l\in R$, $l\in \mathbb{N}$ such that $I_0=R$ and $I_lI_j\subset I_{l+j}\mbox{,\; } \forall l,j\in \mathbb{N}$. Here, $W$ is just a variable to remind us the degree of the ideals $I_l$. Since $\mathcal{G}$ is finitely generated, there exist some $f_1,\ldots ,f_r\in R$ and positive integers (weights) $n_1,\ldots ,n_r\in \mathbb{N}$ such that
	\begin{equation}\label{def:Rees_alg_generadores}
	\mathcal{G}=R[f_1W^{n_1},\ldots ,f_rW^{n_r}]\mbox{.}
	\end{equation}
Rees algebras can be defined over Noetherian schemes $Z$ in the obvious manner.
In this case $\mathcal{G}$ is a sheaf of graded algebras and $I_l$ is a sheaf of ideals for any $l$.
\end{Def}

\begin{Rem}\label{casi_anillos}
Note that this definition is more general than the (usual) of Rees ring where one considers algebras of the form $R[IW]$ for some ideal $I\subset R$. There is another special type of Rees algebras that will play a role in our arguments: these are Rees algebras of the form $R[IW^b]$, for some ideal $I\subset R$ and some positive integer $b\geq 1$. We refer to them as {\em almost Rees rings}.
In fact, every Rees algebra will be \emph{equivalent} in some sense to an almost Rees ring (see Remark \ref{todas_casi} below).
\end{Rem}

\begin{Def} 
Two Rees algebras over a Noetherian ring $R$ are \textit{integrally equivalent} if their integral closure in $\mathrm{Quot}(R)[W]$ coincide. We say that a Rees algebra over $R$, $\mathcal{G}=\oplus _{l\geq 0}I_lW^l$ is \textit{integrally closed} if it is integrally closed as a  ring in $\mathrm{Quot}(R)[W]$. We denote by $\overline{\mathcal{G}}$ the integral closure of $\mathcal{G}$.
\end{Def}

\begin{Rem}\label{todas_casi}
Note that $\overline{\mathcal{G}}$ is also a Rees algebra over $R$, when $R$ is excellent (\cite[\S 1.1]{Br_G-E_V}).
It can be shown that any Rees algebra $\G=\oplus_lI_lW^l$ is finite over an almost Rees ring.
In fact there is some positive integer $N$ such that $\G$ is finite over $R[I_NW^N]$.
Moreover, if $\G$ is finite over $R[I_NW^N]$ then it can be checked that $\G$ is finite over $R[I_LW^L]$ for any $L$ multiple of $N$ (see \cite[Remark 1.3 and Lemma 1.7]{E_V}). 
\end{Rem} 

\begin{Parrafo}
	\textbf{The Singular Locus of a Rees Algebra.} (\cite[Proposition 1.4]{E_V}). Now let $\mathcal{G}$ be a Rees algebra over a smooth scheme $V$ defined over a perfect field $k$. In such case, we can attach a closed set to $\G$, its \textit{singular locus}, Sing$(\mathcal{G})$, by considering all the points $\xi \in V$ such that $\nu _{\xi }(I_l)\geq l$, $\forall l\in \mathbb{N}$. Here $\nu _{\xi}(I)$ denotes the order of the ideal $I$ in the regular local ring $\mathcal{O}_{V,\xi }$. If $\mathcal{G}=R[f_1W^{n_1},\ldots ,f_rW^{n_r}]$, then it can be checked that: 
	$$\mathrm{Sing}(\mathcal{G})=\left\{ \xi \in \mathrm{Spec}(R):\, \nu _{\xi }(f_i)\geq n_i,\; \forall i=1,\ldots ,r\right\} \subset V\mbox{.}$$
\end{Parrafo}

\begin{Ex} \label{Ex:HiperRees} Suppose that $R$ is smooth over a perfect field $k$. 
	Let $X\subset\Spec(R)=V$ be a hypersurface with $I(X)=(f)$ and let $b>1$ be the maximum value of the multiplicity of $X$.
	If we set $\mathcal{G}=R[fW^b]$ then $\Sing(\mathcal{G})=F_b(X)$. Along this paper we will be using a generalization of this description of the maximum multiplicity locus in the case where $X$ is an equidimensional singular algebraic variety (defined over a perfect field $k$) (see Theorem \ref{Th:PresFinita} and the discussion in \ref{PresFinita}). 
\end{Ex}

\begin{Parrafo}
	{\bf Singular locus, integral closure and differential saturation.} 
	A Rees algebra $\mathcal{G}=\oplus _{l\geq 0}I_lW^l$ defined on a smooth scheme $V$ over a perfect field $k$, is \textit{differentially closed} (or {\em differentially saturated}) if there is an affine open covering $\{U_i\}_{i\in I}$ of $V$, such that for every $D\in \mathrm{Diff}^{r}(U_i)$ and $h\in I_l(U_i)$, we have $D(h)\in I_{l-r}(U_i)$ whenever $l\geq r$ (where $\mathrm{Diff}^{r}(U_i)$ is the locally free sheaf over $U_i$ of $k$-linear differential operators of order less than or equal to $r$). In particular, $I_{l+1}\subset I_l$ for $l\geq 0$. We denote by $\mathrm{Diff}(\mathcal{G})$ the smallest differential Rees algebra containing $\mathcal{G}$ (its \textit{differential closure}). (See \cite[Theorem 3.4]{V07} for the existence and construction.)
	
	It can be shown (see \cite[Proposition 4.4 (1), (3)]{V3}) that for a given Rees algebra $\mathcal{G}$ on $V$, 
	$$\Sing(\G) =\Sing (\overline{\G})= \Sing (\mathrm{Diff}(\mathcal{G})).$$
	
\end{Parrafo}
The problem of {\em simplification of the multiplicity of an algebraic variety} can be translated into the problem of {\em resolution of a suitably defined Rees algebra} (see Theorem \ref{Th:PresFinita}). This motivates Definitions \ref{def:transf_law} and \ref{def:res_RA} below (see also Example \ref{ejemplo_hipersuperficie}). 

\begin{Def}\label{def:transf_law}
	Let $\mathcal{G}$ be a Rees algebra on a smooth scheme $V$. A \textit{$\mathcal{G}$-permissible blow up}
$V\stackrel{\pi }{\leftarrow} V_1\mbox{,}$ 
	is the blow up of $V$ at a smooth closed subset $Y\subset V$ contained in $\mathrm{Sing}(\mathcal{G})$ (a {\em permissible center for $\mathcal{G}$}). We denote then by $\mathcal{G}_1$ the (weighted) transform of $\mathcal{G}$ by $\pi $, which is defined as
	$$\mathcal{G}_1:=\bigoplus _{l\in \mathbb{N}}I_{l,1}W^l\mbox{,}$$
	where 
	\begin{equation}\label{eq:transf_law}
	I_{l,1}=I_l\mathcal{O}_{V_1}\cdot I(E)^{-l}
	\end{equation}
	for $l\in \mathbb{N}$ and $E$ the exceptional divisor of the blow up $V\stackrel{\pi }{\leftarrow} V_1$.
\end{Def}

\begin{Def}\label{def:res_RA}
	Let $\mathcal{G}$ be a Rees algebra over a smooth scheme $V$. A \textit{resolution of $\mathcal{G}$} is a finite sequence of transformations
	\begin{equation}\label{diag:res_Rees_algebra}
	\xymatrix@R=0pt@C=30pt{
		V=V_0 & V_1 \ar[l]_>>>>>{\pi _1} & \ldots \ar[l]_{\pi _2} & V_l \ar[l]_{\pi _l}\\
		\mathcal{G}=\mathcal{G}_0 & \mathcal{G}_1 \ar[l] & \ldots \ar[l] & \mathcal{G}_l \ar[l]
	}\end{equation}
	at permissible centers $Y_i\subset \text{Sing} ({\mathcal G}_i)$, $i=0,\ldots, l-1$, such that $\mathrm{Sing}(\mathcal{G}_l)=\emptyset$, and such that the exceptional divisor of the composition $V_0\longleftarrow V_l$ is a union of hypersurfaces with normal crossings. 
\end{Def}

\begin{Ex}\label{ejemplo_hipersuperficie}
	With the setting of Example~\ref{Ex:HiperRees}, a resolution of the Rees algebra $\mathcal{G}=R[fW^b]$ induces a sequence of transformations such that the multiplicity of the strict transform of $X$ descreases:
	\begin{gather*}
	\xymatrix@R=0pt@C=30pt{
		\mathcal{G}=\mathcal{G}_0 & \mathcal{G}_1 \ar[l] & \ldots \ar[l] & \mathcal{G}_{l-1} \ar[l] & \mathcal{G}_l \ar[l]\\
		V=V_0 & V_1 \ar[l]_>>>>>{\pi _1} & \ldots \ar[l]_{\pi _2} & V_{l-1}\ar[l]_{\pi_{l-1}} & V_l \ar[l]_{\pi _l}\\
		\ \; \; \; \; \; \cup & \cup & & \cup & \cup \\
		X=X_0 & X_1 \ar[l]_>>>>>{\pi _1} & \ldots \ar[l]_{\pi _2} & X_{l-1} \ar[l]_{\pi_{l-1}} & X_l \ar[l]_{\pi _l}
	} \\
	b=\max\mult(X_0) = \max\mult(X_1) =\cdots = \max\mult(X_{l-1})>\max\mult(X_l).
	\end{gather*}
	Here each $X_i$ is the strict transform of $X_{i-1}$ after the blow up $\pi_i$.
	Note that the set of points of $X_l$ having multiplicity $b$ is $\Sing(\mathcal{G}_{l})=\emptyset$.
\end{Ex}

\begin{Rem}
Resolution of Rees algebras is known to exists when $V$ is a smooth scheme defined over a field of characteristic zero (\cite{Hir}, \cite{Hir1}). In \cite{V1} and \cite{B-M} different algorithms of resolution of Rees algebras are presented (see also \cite{E_V97}, \cite{E_Hau}). 
\end{Rem}

\begin{Parrafo}\label{Def:HirOrd}
\textbf{Hironaka's order function for Rees algebras.}
(\cite[Proposition 6.4.1]{E_V})
We define the \textit{order of the Rees algebra $\mathcal{G}$ at $\xi \in \mathrm{Sing}(\mathcal{G})$} as:
$$\mathrm{ord}_{\xi }(\mathcal{G}):=\inf _{l\geq 0}\left\{ \frac{\nu_{\xi }(I_l)}{l}\right\} \mbox{.}$$
This is what we call {\em Hironaka's order function of $\G$ at the point $\xi$}. 
If $\mathcal{G}=R[f_1W^{n_1},\ldots ,f_rW^{n_r}]$ and $\xi\in \mathrm{Sing}(\mathcal{G})$ then it can be shown that $\mathrm{ord}_{\xi }(\mathcal{G})=\min _{i=1,\ldots,r}\left\{ \mathrm{ord}_{\xi }(f_iW^{n_i})\right\}$, where $\mathrm{ord}_{\xi}(f_iW^{n_i}):=\frac{\nu_{\xi}(f_i)}{n_i}$, (see \cite[Proposition 6.4.1]{E_V}). Finally, it can be proven that for any point $\xi\in\Sing(\mathcal{G})$ we have
$\ord_{\xi}(\mathcal{G})=\ord_{\xi}(\overline{\mathcal{G}})=\ord_{\xi}(\Diff(\mathcal{G}))$ 
(see \cite[Remark 3.5, Proposition 6.4 (2)]{E_V}).

\

Along this paper we use `$\nu$' to denote the usual order of an element or an ideal at a regular local ring, and `$\ord$' for the order of a Rees algebra at a regular local ring. 
\end{Parrafo}

\

\noindent {\bf Local presentations of the Multiplicity} 
\medskip

\noindent In the following paragraphs we will see that the constructions of Examples \ref{Ex:HiperRees} and 
\ref{ejemplo_hipersuperficie} can be extended to the case in which $X$ is not necesarily a hypersurface. To be more precise, in \cite{V} it is proven that for each $\xi\in F_m(X)$ there is an (\'etale) neighborhood $U\subset X$ of $\xi$ which we denote again by $X$ to ease the notation, and an embedding 
$X\subset V=\Spec(R)$ for some smooth $k$-algebra $R$, together with an $R$-Rees algebra, $\mathcal{G}$, 
so that 
\begin{equation} \label{eq:LocPresSing}
F_m(X)=\Sing(\mathcal{G}). 
\end{equation}
and so that, in addition, given a sequence of blow ups at regular equimultiple centers, 
\begin{equation} \label{eq:LocalGSeq}
\xymatrix@R=0pt@C=30pt{
	V=V_0 & V_1 \ar[l]_>>>>>{\pi _1} & \ldots \ar[l]_{\pi _2} & V_l \ar[l]_{\pi _l}\\
	\ \; \; \; \; \; \cup & \cup & & \cup \\
	X=X_0 & X_1 \ar[l] & \ldots \ar[l] & X_l \ar[l] \\
	\mathcal{G}=\mathcal{G}_0 & \mathcal{G}_1 & \ldots & \mathcal{G}_l
}\end{equation}
 the following equality of closed subsets holds: 
\begin{equation} \label{eq:LocPresSing2}
F_m(X_j)=
%\mathrm{\underline{Max}}\mult_{X_{j}}=
\Sing(\mathcal{G}_j), 
\qquad j=0,1,\ldots,l.
\end{equation}

It is worth mentioning that in fact, the link between $F_m(X)$ and  $\G$ is much stronger (it can be checked that equality (\ref{eq:LocPresSing2}) is also preserved after considering local transformations as in (\ref{local_transformations}).
Thus the problem of finding a simplification of the multiplicity of an algebraic variety is translated into the problem of finding a resolution of a suitable Rees algebra defined on a smooth scheme. 
The local embedding together with the Rees algebra $\G$ strongly linked to $F_m(X)$ is what we call a {\em local presentation of the multiplicity on $X$, $\mult_X$}, which we denote by $(V,\G)$. Precise statements about local presentations can be found for instance in \cite[Part II]{Br_V1} or in \cite{Su_V}.

\begin{Thm} \label{Th:PresFinita}
	\cite[\S 7.1]{V}
	Let $X$ be a reduced equidimensional scheme of finite type over a perfect field $k$.
	Then for every point $\xi\in X$ there exists a local presentation for the function $\mult_X$ in an (\'etale) neighborhood of $\xi$. 
\end{Thm}

We give some ideas about the proof of Theorem \ref{Th:PresFinita} since we will use them in the proof of Theorem \ref{MainTheorem}. 

\begin{Parrafo} \label{PresFinita}
	\textbf{Some ideas behind the proof of Theorem \ref{Th:PresFinita}.\cite[\S5, \S7]{V}} The statement of the theorem is of local nature. So, let us assume that $X$ is an affine algebraic variety of dimension $d$, and let $\xi\in F_m(X)$. Then it can be shown that, after considering a suitably \'etale extension of $B$, which we denote by $B$ again for simplicity, we are in the following setting: 
$X=\Spec(B)$, there is a smooth $k$-algebra, $S$, and a finite extension $S\subset B$ of generic rank $m$, inducing a finite morphism $\delta: \text{Spec}(B)\to \text{Spec}(S)$.
Under these assumptions, $B=S[\theta_1,\ldots, \theta_{n}]$, for some $\theta_1,\ldots, \theta_{n}\in B$ and some $n\in {\mathbb N}$.
Observe that the previous extension induces a natural embedding $X\subset V^{(n+d)}:=\text{Spec}(R)$, where $R=S[x_1,\ldots, x_{n}]$. Let $K(S)$ be the field of fractions of $S$ and let $\mathrm{Quot}(B)$ be the total quotient ring of $B$.
Now, if $f_i(x_i)\in K(S)[x_i]$ denotes the minimal polynomial of $\theta_i$ for $i=1,\ldots,n$, then it can be shown that in fact $f_i(x_i)\in S[x_i]$ and as a consequence $\langle f_1(x_1), \ldots, f_{n}(x_{n})\rangle\subset {\mathcal I}(X)$, where ${\mathcal I}(X)$ is the defining ideal of $X$ in $V^{(n+d)}$.
Finally, if each polinomial $f_i$ is of degree $l_i$, it is proven that the differential Rees algebra
	\begin{equation}
	\label{G_Representa}
	\mathcal{G}^{(n+d)}=\Diff(R[f_1W^{l_1}, \ldots, f_{n}W^{l_{n}}])
	\end{equation}
	is a local presentation of $F_m(X)$ at $\xi$. Observe that the finite morphism $\delta: \text{Spec}(B)\to \text{Spec}(S)$ is transversal with generic rank $m$, the maximum multiplicity of $X$. Therefore conditions (i)-(iii) from Remark \ref{condicion_estrella} hold for all primes in $F_m(X)$. 
\end{Parrafo}	

\begin{Rem} \label{sing_gx}
Local presentations are not unique. For instance, once a local (\'etale) embedding $X\subset V$ is fixed, there may be different ${\mathcal O}_V$-Rees algebras representing $F_m(X)$. 
However, it can be proven that they all lead to the same simplification of the multiplicity of $X$, i.e., they all lead to the same sequence (\ref{eq:LocalGSeq}) with $\Sing\G_l=\emptyset$ (at least in characteristic zero, see \cite{Br_G-E_V}, \cite{Br_V2} and \cite{E_V}).
Moreover, in \cite{Abad} it is proven that the restriction to $X$ of the Rees algebra $\mathcal{G}^{(n+d)}$ defined in (\ref{G_Representa}) is well defined up to integral closure. We denote it by $\G_X$ and refer to it as {\em the ${\mathcal O}_X$-Rees algebra attached to $F_m(X)$}.
Finally, notice that since $\mathcal{G}^{(n+d)}=\oplus J_iW^i$ is a differential Rees algebra, $\Sing(\mathcal{G}^{(n+d)})={\mathbb V}(J_i)$ for all $i\geq 1$ (cf. \cite[Proposition 3.9]{V07}). Therefore, if $\G_X=\oplus I_iW^i$ for suitable ideals $I_i\subset {\mathcal O}_X$, it can be assumed that ${\mathbb V}(I_i)=F_m(X)$ for $i\geq 1$. 
\end{Rem} 

\section{Arcs, jets and Nash multiplicity sequences} \label{Jet_Arcs}

\begin{Def}
	Let $Z$ be a scheme over a field $k$, and let $K\supset k$ be a field extension.
	An {\em m-jet in $Z$} is a morphism $\vartheta: \text{Spec}\left(K[[t]]/\langle t^{m+1}\rangle\right)\to Z$ for some $m\in {\mathbb N}$. 
\end{Def}
If ${\mathcal S}ch/k$ denotes the category of $k$-schemes and ${\mathcal S}et$ the category of sets, then the contravariant functor:
$$\begin{array}{rcl} 
{\mathcal S}ch/k & \longrightarrow & {\mathcal S}et\\
Y & \mapsto & \text{Hom}_k(Y\times_{\text{Spec}(k)}\text{Spec}(k[[t]]/\langle t^{m+1}\rangle), Z)
\end{array}
$$
is representable by a $k$-scheme ${\mathcal L}_m(Z)$, the {\em space of $m$-jets} over $Z$. If $Z$ is of finite type over $k$, then so is ${\mathcal L}_m(Z)$ (see \cite{Vojta}). For each pair $m\geq m'$ there is the (natural) truncation map ${\mathcal L}_{m}(Z) \to {\mathcal L}_{m'}(Z)$. In particular, for $m'=0$, ${\mathcal L}_{m'}(Z)=Z$ and we will denote by ${\mathcal L}_m(Z, \xi)$ the fiber of the (natural) truncation map over a point $\xi\in Z$. Finally, if $Z$ is smooth over $k$ then ${\mathcal L}_m(Z)$ is also smooth over $k$ (see \cite{I2}).

By taking the inverse limit of the ${\mathcal L}_m(Z)$, the {\em arc space of $Z$} is defined, 
$${\mathcal L}(Z):=\lim_{\leftarrow}{\mathcal L}_m(Z).$$ This is the scheme representing the functor (see \cite{Bhatt}): 
$$\begin{array}{rcl} 
{\mathcal S}ch/k & \longrightarrow & {\mathcal S}et\\
Y & \mapsto & \text{Hom}_k(Y\tilde{\times}\text{Spf}(k[[t]]), Z).
\end{array}
$$

A $K$-point in $\L(Z)$ is an {\em arc of $Z$} and can be seen as a morphism $\varphi: \text{Spec}(K[[t]])\to Z$ for some $K\supset k$. The image by $\varphi$ of the closed point is called the {\em center of the arc $\varphi$}. If the center of $\varphi$ is $\xi\in Z$ then it induces a $k$-homomorphism 
${\mathcal O}_{Z,\xi}\to K[t]]$ which we will denote by $\varphi$ too; in this case the image by $\varphi$ of the maximal ideal , $\varphi({\mathfrak m}_{\xi})$, generates an ideal $\langle t^l\rangle\subset K[[t]]$ and then we will say that {\em the order of $\varphi$ is $l$} and we will denote it by $\nu_t(\varphi)$. We will denote by $\L(Z, {\xi})$ the set of arcs in $\L(Z)$ with center $\xi$. The {\em generic point of $\varphi$ in $Z$} is the point in $Z$ determined by the $\text{kernel}$ of $\varphi$.

\begin{Def}
An arc $\varphi: \text{Spec}(K[[t]])\to Z$ is {\em thin} if it factors through a proper closed suscheme of $Z$.
Otherwise we say that $\varphi$ is {\em fat}.
\end{Def}

\

\noindent {\bf Nash multiplicity sequences}
\medskip

\noindent Let $X$ be an algebraic variety defined over a perfect field $k$ and let $\xi \in X$ be a (closed) point.
Assume that $X$ is locally a hypersurface in a neighborhood of $\xi$, $X\subset V$, where $V$ is smooth over $k$, and work at the completion $\widehat{\mathcal O}_{V,{\mathfrak m}_{\xi}}$. Under these hypotheses, in \cite{L-J}, Lejeune-Jalabert introduced the {\em Nash multiplicity sequence along an arc $\varphi\in {\mathcal L}(X, {\xi})$} (in fact, the hypotheses in \cite{L-J} are weaker, but we are interested in working over perfect fields). The Nash multiplicity sequence along $\varphi$ is a non-increasing sequence of non-negative integers 
\begin{equation}
\label{Nash_sequence}
m_0\geq m_1\geq \ldots \geq m_l=m_{l+1}=...\geq 1, 
\end{equation}
where $m_0$ is the usual multiplicity of $X$ at $\xi$, and the rest of the terms are computed by considering suitable stratifications on ${\mathcal L}_{m}(X, {\xi})$ defined via the action of certain differential operators on the fiber of the jets spaces ${\mathcal L}_{m}(\text{Spec}(\widehat{\mathcal O}_{V,{\mathfrak m}_{\xi}}))$ over $\xi$ for $m\in {\mathbb N}$. The sequence (\ref{Nash_sequence}) can be interpreted as the {\em multiplicity of $X$ along the arc $\varphi$}: thus it can be seen as a refinement of the usual multiplicity. The sequence stabilizes at the value given by the multiplicity $m_l$ of $X$ at the generic point of the arc $\varphi $ in $X$ (see \cite[\S 2, Theorem 5]{L-J}). 

\vspace{0.2cm}

In \cite{Hickel93}, Hickel generalized Lejeune's construction to the case of an arbitrary variety $X$,  and in \cite{Hickel05} he presented the sequence (\ref{Nash_sequence}) in a different way which we will explain along the following lines.

\vspace{0.2cm}

Since the arguments are of local nature, let us suppose that $X=\text{Spec}(B)$ is affine.
Let $\xi\in X$ be a point (which we may assume to be closed) of multiplicity $m$, and let $\varphi: B\to K[[t]]$ be an arc in $X$ centered at $\xi$. Consider the natural morphism
$$\Gamma _0=\varphi \otimes i:B\otimes_k k[t]\rightarrow K[[t]]\mbox{,}$$
which is additionally an arc in $X_0=X\times \mathbb{A}^1_k$ centered at the point $\xi _0=(\xi ,0)\in X_0$. This arc determines a sequence of blow ups at points:
\begin{equation}\label{intro:diag:Nms}
\xymatrix@R=15pt@C=50pt{
	\mathrm{Spec}(K[[t]]) \ar[dd]^{\Gamma _0} \ar[ddr]^{\Gamma _1} \ar[ddrrr]^{\Gamma _l} & & & & \\
	& & & & \\
	X_0=X\times \mathbb{A}^1_k & X_1 \ar[l]_>>>>>{\pi _1} & \ldots \ar[l]_{\pi _2} & X_l \ar[l]_{\pi _l} & \ldots \\
	\xi _0=(\xi ,0) & \xi _1 & \ldots & \xi _l & \ldots 
}\end{equation}
Here, $\pi _i $ is the blow up of $X_{i-1}$ at $\xi _{i-1}$, where $\xi _{i}=\mathrm{Im}(\Gamma _{i})\cap \pi _{i}^{-1}(\xi _{i-1})$ for $i=1,\ldots ,l,\ldots $, and $\Gamma _i$ is the (unique) arc in $X_i$ with center $\xi _i$ which is obtained by lifting $\Gamma _0$ via the proper birational morphism $\pi _i\circ \ldots \circ \pi _1$. This sequence of blow ups defines a non-increasing sequence
\begin{equation}
\label{Nash_sequence_2}
m=m_0\geq m_1\geq \ldots \geq m_l=m_{l+1}=...\geq 1, 
\end{equation}
where $m_i$ corresponds to the multiplicity of $X_i$ at $\xi _i$ for each $i=0,\ldots ,l,\ldots $. Note that $m_0$ is nothing but the multiplicity of $X$ at $\xi $, and it is proven that for hypersurfaces the sequence (\ref{Nash_sequence_2}) coincides with the sequence (\ref{Nash_sequence}) above. We will refer to the sequence of blow ups in (\ref{intro:diag:Nms}) as the {\em sequence of blow ups directed by $\varphi$}.

\

\noindent {\bf The persistance}
\medskip

\begin{Def}\label{def:rho}
	Let $\varphi $ be an arc in $X$ with center $\xi \in X$, a point of multiplicity $m>1$. Suppose that the 
	generic point of $\varphi$ is not contained in the stratum of points of multiplicity $m$ of $X$. We denote by $\rho _{X,\varphi }$ the minimum number of blow ups directed by $\varphi $ which are needed to lower the Nash multiplicity of $X$ at $\xi $. That is, $\rho _{X,\varphi }$ is such that $m=m_0=\ldots =m_{\rho _{X,\varphi }-1}>m_{\rho _{X,\varphi }}$ in the sequence (\ref{Nash_sequence_2}) above. We call $\rho _{X,\varphi }$ the \textit{persistance of $\varphi$} (we will see in Remark \ref{persistencia_g} that the persistance is always finite). 
\end{Def}

\begin{Rem}\label{Hickel_construction}
	Using Hickel's construction, it can be checked that the first index $i\in \{1,\ldots, l+1\}$ for which there is a strict inequality in (\ref{Nash_sequence_2}) (i.e., the first index $i$ for which $m_0>m_i$) can be interpreted as the minimum number of blow ups needed to {\em separate the graph of $\varphi$ from } the stratum of points of multiplicity $m_0$ of $X_0$ (actually, to be precise, this statement has to be interpreted in $B\otimes K[[t]]$, where the graph of $\varphi$ is defined). 
\end{Rem}

Next we define a normalized version of $\rho _{X,\varphi }$ in order to avoid the influence of the order of the arc in the number of blow ups needed to lower the Nash multiplicity. 

\begin{Def}\label{def:rho_bar} 
	For a given arc $\varphi: \text{Spec}(K[[t]])\to X$ with center $\xi\in X$, we will write
	$$\bar{\rho }_{X,\varphi }=\frac{\rho _{X,\varphi }}{\nu_t(\varphi )}\mbox{.} $$
\end{Def}

\begin{Def} \label{rho_function} For each point $\xi\in X$ we define the functions: 
	\begin{equation}
	\label{funcion_rho_1}
	\begin{array}{rrclcrrcl}
	\rho_X: & \L(X, \xi)& \to & {\mathbb Q}_{\geq 0} \cup \{\infty\} & \ \ \text{ and } \ \ & \overline{\rho}_X: & \L(X, \xi) & \to & {\mathbb Q}_{\geq 0 } \cup \{\infty\}\\ 
	& \varphi & \mapsto & \rho_{X,\varphi} & & 
	& \varphi & \mapsto & \overline{\rho}_{X,\varphi}.
	\end{array} 
	\end{equation} 
\end{Def}

\begin{Rem} \label{rho_etale} Many of our arguments will be developed, locally, in an \'etale neighborhood of a point $\xi\in X$, but the persistance is stable after considering \'etale morphims. 
	In fact the whole sequence $\{m_i\}_{i\geq 0}$ in (\ref{Nash_sequence_2}) does not change in an \'etale neighborhood of $\xi\in X$ in the following sense. Suppose $\mu: \widetilde{X} \to X$ is an \'etale morphism with $\mu(\widetilde{\xi})=\xi$, and let $\varphi: \Spec(K[[t]]) \to X$ be an arc with center $\xi$. Then there is a lifting with center $\widetilde{\xi}$, $\widetilde{\varphi}: \Spec(\widetilde{K}[[t]]) \to \widetilde{X}$, where $\widetilde{K}$ is a separable extension of $K$. 
	If the Nash multiplicity sequence for the arc $\widetilde{\varphi}$ is $\{\widetilde{m}_i\}_{i\geq 0}$, and the Nash multiplicity sequence for $\varphi$ is $\{m_i\}_{i\geq 0}$, then it can be checked that $m_i=\widetilde{m}_i$ for all $i\geq 0$.
In particular the persistance of $\varphi$ is the same as the persistance of $\widetilde{\varphi}$, and so does the normalized persistance at $\varphi$ and $\widetilde{\varphi}$, i.e., 
$\rho_{X,\varphi}=\rho_{\widetilde{X},\widetilde{\varphi}} \ \text{ and } \ \overline{\rho}_{X,\varphi}=\overline{\rho}_{\widetilde{X},\widetilde{\varphi}}$. 
We refer to \cite[Remark 2.8]{BEP3} for full details. 
\end{Rem}

\

\noindent{\bf The ${\mathbb Q}$-persistance}
\medskip

\begin{Def}\label{q_persistencia}
Let $\varphi $ be an arc in $X$ with center $\xi \in X$, a point of multiplicity $m>1$, $\varphi :\mathrm{Spec}(K[[t]])\longrightarrow X$. Consider the family of arcs given as $\varphi _n=\varphi \circ i_n$ for $n>1$, where $i_n^*:K[[t]]\longrightarrow K[[t^n]]$ is the $K$-morphism that maps $t$ to $t^n$. Then the {\em ${\mathbb Q}$-persistance of $\varphi$}, ${r}_{X,\varphi}$, is defined as the limit: 
	\begin{equation}\label{r_limite_n} 
	{r}_{X,\varphi }:= \lim _{n\rightarrow \infty }\frac{\rho _{X,\varphi _{n}}}{n}. 
	\end{equation}
	And the {\em normalized ${\mathbb Q}$-persistance of $\varphi$} is:
	\begin{equation}\label{r_limite} 
	\bar{r}_{X,\varphi }:=\frac{r_{X,\varphi}}{\nu_t({\varphi})}=\frac{1}{\nu_t(\varphi )}\cdot \lim _{n\rightarrow \infty }\frac{\rho _{X,\varphi _{n}}}{n}. 
	\end{equation}	
\end{Def}

 As we will see in Remark \ref{persistencia_g} below, the ${\mathbb Q}$-persistance of $\varphi$ can be interpreted as the {\em order of contact of the arc $\varphi$ with the stratum of multiplicity $m_0$ of the variety $X_0$} (see expression \ref{r_escrito}). There we will also justify that both limits (\ref{r_limite_n}) and (\ref{r_limite}) exist.

\begin{Def} \label{r_function}
For each point $\xi\in X$ we define the functions: 
\begin{equation}
\label{funcion_r}
\begin{array}{rrclcrrcl}
	r_X: & \L(X, \xi)& \to & {\mathbb Q}_{\geq 0} \cup \{\infty\} & \ \ \text{ and } \ \ & \overline{r}_X: & \L(X, \xi) & \to & {\mathbb Q}_{\geq 0 } \cup \{\infty\}\\ 
	& \varphi & \mapsto & r_{X,\varphi} & & 
	& \varphi & \mapsto & \overline{r}_{X,\varphi}.
	\end{array} 
\end{equation}
\end{Def}

\begin{Rem} \label{persistencia_g} 
Let $\varphi\in \L(X, {\xi})$ and suppose that $\G_X$ is defined on $X$. Then, it can be shown that: 	
\begin{equation}
\label{r_escrito}
	r_{X,\varphi}=\mathrm{ord}_{t}(\varphi (\G_{X}))\in \mathbb{Q}_{\geq 1}, 
\end{equation}
and hence, 
\begin{equation}
\label{r_escrito_n}
	\bar{r}_{X,\varphi}=\frac{\mathrm{ord}_{t}(\varphi (\G_{X}))}{\nu_t(\varphi)}\in \mathbb{Q}_{\geq 1}\mbox{,}
\end{equation}
where, if we assume that $\G_{X}$ is generated by $g_1W^{b_1},\ldots , g_rW^{b_r}$ in some affine chart $\Spec(B)$ of $X$ containing the center of the arc $\varphi: B\to K[[t]]$, then
\begin{equation}
\label{arco_algebra} 
	\varphi({\mathcal{G}}_{X}):=K[[t]][\varphi (g_1)W^{b_1},\ldots ,\varphi(g_{r})W^{b_{r}}]\subset K[[t]][W].	
\end{equation} 
See \cite[Corollary 4.3.4]{P-ET}, and \cite[Proposition 5.9, Remark 5.10 and \S 6]{BEP2}. 
From here it can be checked that, if the generic point of the arc $\varphi$ is not contained in $F_m(X)=\Sing (\Gn)$, then $\varphi(\G_{X})\subset K[[t]]$ is a non zero Rees algebra. As a consequence, $r_{X,\varphi}$ is finite. Moreover, it can be shown that $\rho _{X,\varphi }$ is obtained by taking the integral part of $r_{X,\varphi }$ (see \cite[Proposition 5.11]{BEP2}, and also \cite{Br_E_P-E}), and therefore it is also a finite number. 
	
Now, notice that the expression (\ref{r_escrito}) can be computed in an \'etale neighborhood $\widetilde{X}$ of $\xi\in X$ where $\G_X$ is defined (see \ref{PresFinita}). If $\varphi\in \L(X, {\xi})$ then there is always a lifting $\widetilde{\varphi} \in \L(\widetilde{X}, \widetilde{\xi})$ as in Remark \ref{rho_etale} with the same Nash multiplicity sequence. Hence, 
\begin{equation}
\label{r_etale}
\bar{r}_{X,\varphi}=\frac{1}{\nu_t(\varphi )}\cdot \lim _{n\rightarrow \infty }\frac{\rho _{X,\varphi _{n}}}{n}=\frac{1}{\nu_t(\widetilde{\varphi})}\cdot \lim _{n\rightarrow \infty }\frac{\rho _{X,\widetilde{\varphi} _{n}}}{n}=\frac{\mathrm{ord}_{t}(\widetilde{\varphi} (\G_{X}))}{\nu_t(\widetilde{\varphi})}=\bar{r}_{\widetilde{X},\widetilde{\varphi}}.
\end{equation} 	
From here it also follows that 	the functions $\rho_{X,}$ and $r_{X}$ from Definitions \ref{rho_function} and \ref{r_function} encode the same information. 
\end{Rem} 

\begin{Parrafo} \label{entero_arcos}{\bf Integral closure of Rees algebras and arcs.}
Let $k$ be a field, let $B$ be a (not necessarily smooth) reduced  excellent  $k$-algebra, and let $\G$ be a Rees algebra over $B$.
Set $X=\text{Spec}(B)$.
For any arc $\varphi\in {\mathcal L}(X)$, $\varphi: B\to K[|t|]$, with $k\subset K$ a extension field, the image via $\varphi$ of $\G$ generates a Rees algebra over $K[|t]]$.
It can be checked (see \cite[4.6]{BEP2}) that for all arcs $\varphi\in {\mathcal L}(X)$, 
\begin{equation}
\label{igualdad}
\ord_t(\varphi(\G))=\ord_t(\varphi(\overline{\G})).
\end{equation}

On the other hand, given two Rees algebras $\G$ and $\G'$ on $X$, it can be shown that if for any fat arc $\varphi\in \LX$, 
$\ord_t(\varphi(\G))=\ord_t(\varphi({\G'}))$, then $\overline{\G}=\overline{\G'}$. This follows from the fact that there are ideals $I, J\subset {\mathcal O}_X$ such that, up to integral closure it can be assumed that $\G={\mathcal O}_X[IW^b]$ and $\G={\mathcal O}_X[JW^b]$ for some positive integer $b$ (see Remark \ref{todas_casi}). Thus $\overline{\G}=\overline{\G'}$ if an only if $\overline{I}=\overline{J}$. Now our hypothesis implies that $\nu_t(\varphi({\G}))/b=\nu_t(\varphi({\G'}))/b$ for all fat arcs $\varphi\in \LX$. And now the claim follows from the valuative criterion for integral closure of ideals.
 
\end{Parrafo}

\section{Proof of Theorem \ref{MainTheorem}} \label{proof}

\begin{proof}
The statement is of local nature so we may assume that $X=\text{Spec}(B)$ and $X'= \text{Spec}(B')$ are an affine algebraic varieties   of dimension $d$. Let $\xi'\in F_{rm}(X')$, and let $\xi=\beta(\xi')$. Then $\xi\in F_m(X)$. Arguing as in \ref{PresFinita}, after considering a suitably defined \'etale extension of $B$, which we denote by $B$ again, we may assume that there is a finite morphism $\delta: X \to \Spec(S)=V^{(d)}$, of generic rank $m$, with $S$ a smooth $k$-algebra of dimension $d$, an immersion $X \hookrightarrow V^{(n+d)}=\Spec(S[x_1,\ldots, x_n])$ and an $\mathcal{O}_{V^{(n+d)}}$-Rees algebra $\G^{(n+d)}$ (see (\ref{G_Representa})), which we assume to be differentially closed, representing the multiplicity of $X$. The \'etale extension of $B$ induces an \'etale extension of $B'$ which we denote by $B'$ too. Recall that transversality is preserved by \'etale morphisms (see Remark \ref{etale_transversal}). Thus, since $X'\to X$ is transversal, it follows that the induced extension $S\subset B'$ is transversal too. Hence we have the following diagram:
\begin{equation}\label{diagrama_presentacion}
\xymatrix@R=20pt@C=20pt{R'=S[x_1,\ldots ,x_{n}, x_{n+1},\ldots, x_{n'}] \ar[r] & A' \ar[rr] & \ & B'=S[\theta_1,\ldots, \theta_n, \theta_{n+1},\ldots, \theta_{n'}] \\
	R=S[x_1,\ldots ,x_{n}] \ar[r] \ar[u] & A \ar[rr] \ar[u] & \ & B= S[\theta_1,\ldots, \theta_n] \ar[u]^{\beta^*}\\
	S \ar[u]^{\alpha^*} 
	\ar[urrr]_{\delta^*} 
	& & & 
}
\end{equation}	
	where $A=S[x_1,\ldots ,x_{n}] /\langle f_1,\ldots, f_{n}\rangle$, $A'=S[x_1,\ldots ,x_{n}, x_{n+1},\ldots, x_{n'}] /\langle f_1,\ldots, f_{n}, f_{n+1}, \ldots, f_{n'} \rangle$ and each $f_i(x_i)\in S[x_i]$ is the minimum polynomial of $\theta_i$ over $K(S)$ for $i=1,\ldots, n, n+1,\ldots, n'$. Therefore,   the differential $R'$-Rees algebra  
\begin{equation}
\label{G_Representa_prima}
\G^{(d+n')}:=\Diff(R'[f_1W^{l_1}, \ldots, f_{n}W^{l_{n}}, f_{n+1}W^{l_{n+1}}, \ldots,f_{n'}W^{l_{n'}}]) 
\end{equation}
and the differential $R$-Rees algebra
\begin{equation}
\label{G_Representa_prima_no}
\G^{(d+n)}:=\Diff(R'[f_1W^{l_1}, \ldots, f_{n}W^{l_{n}}]) 
\end{equation}
 represent the maximum multiplicity of $X'$ and $X$ respectively. Observe that there is a natural inclusion of $R'$-Rees algebras, 
 $ \G^{(d+n)} R'\subset \G^{(d+n')}$, that induces a natural inclusion $\G_X\subset \G_{X'}$. 
 
 Now, by Remark \ref{condicion_estrella}, ${\mathfrak m}_{\delta(\xi)}B$ is a reduction of ${\mathfrak m}_{\xi}$, and ${\mathfrak m}_{\xi}B'$ is a reduction of ${\mathfrak m}_{\xi'}$. 
 Thus, for an arc $\varphi'\in \L(X', {\xi'})$ 	
 \begin{equation}
 \label{mismo_orden} 
 \nu_t(\varphi')=\nu_t(\beta_{\infty}(\varphi')). 
 \end{equation}
	
\noindent ($\Rightarrow$) Suppose that $\G_X\subset \G_{X'}$ is a finite extension of Rees algebras.
Then by \ref{entero_arcos}, for any arc $\varphi'\in \LXP$ with center a point $\xi' \in F_{rm}(X')$, $\ord_t(\beta_{\infty}(\varphi')(\G_X))=\ord_t(\varphi'(\G_{X'}))$.
By (\ref{mismo_orden}), $\nu_t(\varphi'({\mathfrak m_{\xi'}}))
=\nu_t(\beta_{\infty}(\varphi')({\mathfrak m_{\beta(\xi')}}))$. Thus the conclusion follows from (\ref{r_escrito_n}) and (\ref{rho_r}). 

\

\noindent ($\Leftarrow$) Assume now that for each arc $\varphi'\in {\mathcal L}(X')$ with center contained in
$F_{rm}(X')$, the equality 
$$\rho_{X',\varphi'}=\rho_{X,\beta_{\infty}(\varphi')}$$
holds. Then by (\ref{rho_r}), for all arcs with center in $F_{rm}(X')$,
$$r_{X',\varphi'}=r_{X,\beta_{\infty}(\varphi')}.$$
From here it follows that for these arcs, $\ord_t(\beta_{\infty}(\varphi')(\G_X))=\ord_t(\varphi'(\G_{X'}))$.
On the other hand, if the center of an arc $\varphi'\in \LXP$ is not contained in $F_{rm}(X')$, then  by the hypotheses of the theorem, the center of $\beta_{\infty}(\varphi')$ is not in $F_m(X)$,   and
by Remark \ref{sing_gx},  $\ord_t(\varphi'(\G_{X'}))=\ord_t(\beta_{\infty}(\varphi')(\G_X))=0$.  Thus the conclusion follows from \ref{entero_arcos}. 
\end{proof}

\end{document}